\documentclass[oneside,a4paper,12pt,reqno]{amsproc}

\usepackage{amssymb}
\usepackage[english]{babel}
\usepackage[T1]{fontenc}
\usepackage[latin1]{inputenc}
\usepackage{fullpage}

\DeclareMathOperator{\PGL}{PGL}
\DeclareMathOperator{\PSL}{PSL}
\DeclareMathOperator{\PGammaL}{P\Gamma L}

\DeclareMathOperator{\PHG}{PHG}
\DeclareMathOperator{\PG}{PG}
\DeclareMathOperator{\AG}{PG}
\DeclareMathOperator{\GR}{GR}
\DeclareMathOperator{\rad}{rad}

\DeclareMathOperator{\Aut}{Aut}

\theoremstyle{plain}
\newtheorem{thm}{Theorem}[section]

\theoremstyle{definition}

\usepackage{amsfonts,latexsym,graphicx}
\usepackage{hyperref}
\usepackage{times}
\usepackage{enumerate}

\newcommand{\F}{\mathbb{F}}
\newcommand{\N}{\mathbb{N}}
\newcommand{\Z}{\mathbb{Z}}
\newcommand{\G}{\mathbb{G}}
\newcommand{\IS}{\mathbb{S}}
\newcommand{\T}{\mathbb{T}}

\newcommand{\arcone}{\mathfrak{k}_{\Z_{25}}}
\newcommand{\arctwo}{\mathfrak{k}_{\IS_5}}

\providecommand{\abs}[1]{\lvert#1\rvert}

\begin{document}

\title{New Complete $2$-Arcs in the Uniform Projective Hjelmslev Planes Over Chain Rings of Order $25$}

\author{Michael Kiermaier}
\address{Michael Kiermaier\\Mathematical Department\\University of Bayreuth\\D-95440 Bayreuth\\Germany}
\email{michael.kiermaier@uni-bayreuth.de}
\urladdr{http://www.mathe2.uni-bayreuth.de/michaelk/}

\author{Matthias Koch}
\address{Matthias Koch\\Mathematical Department\\University of Bayreuth\\D-95440 Bayreuth\\Germany}
\email{matthias.koch@uni-bayreuth.de}
%\urladdr{}


\begin{abstract}
In this paper a $2$-arc of size $21$ in the projective Hjelmslev plane $\PHG(2,\Z_{25})$ and a $2$-arc of size $22$ in $\PHG(2,\F_5[X]/(X^2))$ are given.
Both arcs are bigger than the $2$-arcs previously known in the respective plane.
Furthermore, we will give some information on the geometrical structure of the arcs.
\end{abstract}

\maketitle

\section{Introduction}
It is well known that a Desarguian projective plane of order $q$ admits a $2$-arc of size $q + 2$ if and only if $q$ is even.
These $2$-arcs are called \emph{hyperovals}.
The biggest $2$-arcs in the Desarguian projective planes of odd order $q$ have size $q + 1$ and are called \emph{ovals}.

For an uniform projective Hjelmslev plane over a chain ring $R$ of size $q^2$, the situation is somewhat similar:
In $\PHG(R,2)$ there exists an hyperoval -- that is a $2$-arc of size $q^2 + q + 1$ -- if and only if $G$ is a Galois ring of size $q^2$ with $q$ even, see \cite{Honold-Landjev-2005-FFA11:292-304, Honold-Landjev-2001-DM231:265-278}

In the remaining uniform projective Hjelmslev planes, the siutation is less clear. It is known \cite{Honold-Landjev-2001-DM231:265-278} that for even $q$ and $R$ not a Galois ring, a $2$-arc has at most size $q^2 + q$, and for odd $q$ a $2$-arc has at most size $q^2$.
For $\#R\leq 16$ the exact values were determined either by theory or computationally \cite{Honold-Landjev-2001-DM231:265-278, Kiermaier-2006, Honold-Kiermaier-2006-ProcACCT10:112-117}, leaving $\#R = 25$ as the smallest case where the exact sizes $n_2(R)$ of a maximum $2$-arc in $\PHG(2,R)$ are not known.
Up to isomorphism there are two finite chain rings of size $25$ of composition length $2$, these are $\Z_{25}$ and $\F_5[X]/(X^2)$.
When we started our search, the biggest known $2$-arcs had size $20$ \cite{Honold-Landjev-2001-DM231:265-278} or $18$ \cite{Boumova-Landjev-2004}, respectively.

In the following two sections, we give a brief introduction to finite chain rings and Hjelmslev geometries.
For details, see for example \cite{Honold-Landjev-2001-DM231:265-278} and the references cited there.

\section{Finite Chain Rings}
A ring\footnote{Rings are assumed to contain an unity element and to be associative, but not necessarily commutative.} $R$ is called \emph{chain ring} if the lattice of the left-ideals is a chain.
A chain ring is necessarily local, so there is a unique maximum ideal $N = \rad(R)$ and the quotient ring $R/N$ is a division ring. In the following we will only consider finite chain rings, where we get $R/N\cong \F_q$ with a prime power $q = p^r$, $p$ prime.
We will need the projection $\phi : R\rightarrow \F_q$, $a\mapsto a\mod N$, which is a surjective ring homomorphism.

The number of ideals of $R$ reduced by $1$ is the composition length of $R$, considered as a left module ${}_R R$. This number will be denoted by $m$.

An important subclass of the finite chain rings are the \emph{Galois rings}. Their definition is a slight generalization of the construction of finite fields via irreducible polynomials:

Let $p$ be prime, $r$ and $m$ positive integers, $q = p^r$ and $f\in\Z_{p^m}[X]$ be a monic polynomial of degree $r$ such that the image of $f$ modulo $p$ is irreducible in $\F_p[X]$.
Then the \emph{Galois ring} of order $q^m$ and characteristic $p^m$ is defined as
\[
    \GR(q^m,p^m) = \Z_{p^m}[X]/(f)
\]
Up to isomorphism, the definition is independent of the exact choice of $f$.
The symbols $p$, $q$, $r$ and $m$ are consistent with the earlier definitions: It holds $\GR(q^m,p^m)/\rad(\GR(q^m,p^m)) \cong \F_q$ and the composition length of $\GR(q^m,p^m)$ is $m$.
Furthermore, the Galois rings contain the finite fields and the integer residues modulo a prime power:
$\GR(p^m,p^m)\cong\Z_{p^m}$ and $\GR(p^r,p)\cong \F_{p^r}$.

While the fields are exactly the chain rings of composition length $1$, in this article we are interested in finite chain rings of composition length $2$. The isomorphism types of these rings are known:

\begin{thm}[see \cite{Cronheim-1978-GeomDed7:287-302}]
Let $R$ be a finite chain ring of composition length $2$, $N = \rad(R)$ and $R/N\cong\F_q$.
Then $\#R = q^2$, and exactly one of the following statements is true:
\begin{enumerate}[(a)]
    \item $R$ is isomorphic to the Galois ring $\GR(q^2,p^2)$ of order $q^2$ and characteristic $p^2$.
    \item There is an unique automorphism $\sigma$ of $\F_q$ such that $R$ is isomorphic to the \emph{$\sigma$-dual numbers} $\F_q[X,\sigma]/(X^2)$.\footnote{$\F_q[X,\sigma]$ is a \emph{skew polynomial ring} over $\F_q$. Its addition is defined as in the usual polynomial ring, and the multiplication is the distributive extension of the rule $X \lambda = \sigma(\lambda) X$ for each scalar $\lambda\in\F_q$.}
\end{enumerate}
\end{thm}

We see that there are $r + 1$ isomorphism classes of chain rings of composition length $2$ and order $q^2$.
Among these rings $2$ are commutative, namely $\GR(q^2,p^2)$ and $\IS_q := \F_q[X]/(X^2)$.
The Galois ring is the unique one with characteristic $p^2$, all the others have characteristic $p$. The smallest such chain rings are:

\[
\begin{array}{r|ccc}
q & \multicolumn{3}{c}{R}\\
\hline
2 & \Z_4 & \IS_2 = \F_2[X]/(X^2) &\\
3 & \Z_9 & \IS_3 = \F_3[X]/(X^2) &\\
4 & \G_4 := \GR(16,4) & \IS_{4} = \F_4[X]/(X^2) & \T_{4} := \F_4[X,a\mapsto a^2]/(X^2)\\
5 & \Z_{25} & \IS_{5} = \F_5[X]/(X^2) &
\end{array}
\]
The only non-commutative ring in this list is $\T_4$.

\section{Arcs in Projective Hjelmslev Planes}                                                                 The \emph{projective Hjelmslev geometry} $\PHG(k,R)$ of dimension $k$ over a finite chain ring $R$ is defined as follows: The point set $\mathcal{P}(\PHG(k,R))$ [line set $\mathcal{L}(\PHG(k,R))$] is the set of the free rank $1$ [rank $k$] right submodules of the module $R^{k+1}$, and the incidence is given by set inclusion.

We have $\abs{\mathcal{P}(\PHG(k,R))} = \abs{\mathcal{L}(\PHG(k,R))} = \frac{q^{k+1} - 1}{q-1} q^{k(m-1)}$.
For $m = 1$, $R$ is a finite field and $\PHG(k,R)$ is the classical projective geometry $\PG(R,k)$ of dimension $k$ over $R$.
For $m > 1$ however, two different lines in $\PHG(k,R)$ may meet in more than one point.

The map $\phi$, extended to $R^{k+1}$, is a collineation $\PHG(k,R)\rightarrow \PG(k,\F_q)$.
Let $P,Q\in\mathcal{P}(\PHG(k,R))$. There is more than one line passing through $P$ and $Q$ if and only if $\phi(P) = \phi(Q)$. The preimages $\phi^{-1}(P)$ [$\phi^{-1}(L)$] with $P\in\mathcal{P}(\PG(k,\F_q))$ [$L\in\mathcal{L}(\PG(k,\F_q))$] are called \emph{point} [\emph{line}] \emph{neighbor classes} of $\PHG(k,R)$.
The restriction of a projective Hjelmslev geometry to a single point neighbor class is isomorphic to the affine geometry $\AG(k,\F_q)$.
The group of collineations of $\PGL(k,R)$ is exactly the semilinear projective group $\PGammaL(k+1,R)$ \cite{Kreuzer-1988}.

In the following, we restrict us to the projective Hjelmslev \emph{planes} $\PHG(2,R)$.
If $R$ has composition length $2$, such a plane is called \emph{uniform}.
For $n\in\N$, a set of points $\mathfrak{k}\subseteq\mathfrak{P}(\PHG(2,R))$ of size $n$ is called \emph{projective $(n,u)$-arc}, if some $u$ elements of $\mathfrak{k}$ are collinear, but no $u+1$ elements of $\mathfrak{k}$ are collinear.
If we allow $\mathfrak{k}$ to be a \emph{multi}set of points in this definition\footnote{Of course we have to
respect multiplicities for counting the number of collinear points.}, $\mathfrak{k}$ is called \emph{$(n,u)$-multiarc}.
We denote by $n_u(R)$ the maximum size of an $u$-multiarc in the projective Hjelsmelv plane $\PHG(2,R)$.

In the case $u = 2$ the discrimination of projective arcs and multiarcs is not important, since the only proper $2$-multiarc is a single point of multiplicity $2$.
So we will simply use the expression $2$-arc.

For a $2$-arc $\mathfrak{k}$ in an uniform Hjelmslev plane over a chain ring with odd parameter $q$, it is known that the complement of the image $\phi(\mathfrak{k})$ is a blocking set in $\PG(2,\F_q)$.

%\cite{Honold-Landjev-2001-DM231:265-278} contains a table for arcs over chain rings of composition length $m=2$ and order $\leq 25$.
%A few new arcs can be found in \cite{Hemme-Honold-Landjev-2000}, and further improvements for chain rings of composition length $m=2$ and order $9$ and $25$ are published in \cite{Boumova-Landjev-2004}.
%In \cite{Kiermaier-2006} a complete classification of $(n,u)$-multiarcs was done for small $u$ in small Hjelmslev geometries over chain rings $\neq \Z_{16}$ of order $\leq 16$, which again yielded some improvements of th
%e bounds.                                                                                                     The most important results of this search can also be found in \cite{Honold-Kiermaier-2006}.

The following table shows the known values of $n_2(R)$ for the finite chain rings $R$ with $m = 2$ and $\#R\leq 25$.
The values for the rings $R$ with $\# R\leq 9$, as well as the lower bound for $\IS_4$ and the upper bounds for $\G_4$, $\Z_{25}$ and $\IS_5$ were given in \cite{Honold-Landjev-2001-DM231:265-278}.
The lower bound for $\G_4$ was given in \cite{Hemme-Honold-Landjev-2000}, and the upper bound for $\IS_4$ as well as lower and upper bound for $\T_4$ can be found in \cite{Honold-Kiermaier-2006-ProcACCT10:112-117}.
For the chain rings $R$ with $q = 5$ the table shows a range, since the exact value of $n_2(R)$ is not known. 
The lower bounds of these ranges are improved by our search.

\[
\begin{array}{c||cc|cc|ccc|cc}
q & \multicolumn{2}{|c|}{2} & \multicolumn{2}{|c|}{3} & \multicolumn{3}{|c|}{4} & \multicolumn{2}{|c}{5}\\
\hline
R & \Z_4 & \IS_2 & \Z_9 & \IS_3 & \G_4 & \IS_4 & \T_4 & \Z_{25} & \IS_5\\
\hline
n_2(R) & 7 & 6 & 9 & 9 & 21 & 18 & 18 & \mathbf{21}-25 & \mathbf{22}-25
\end{array}
\]

\section{The New Arcs}
In this section we give the new arcs and some analysis of their geometrical structure.
The Hjelmslev planes $\PHG(2,\Z_{25})$ and $\PHG(2,\IS_5)$ both consist of $775$ points and lines, and of $31$ point and line neighbor classes. A single neighbor class contains $25$ points respectively lines.

\subsection{A $(21,2)$-arc in $\PHG(2,\Z_{25})$}
A $(21,2)$-arc $\arcone$ in $\PHG(2,\Z_{25})$ is given by the points
\[
\begin{array}{ccc}
(1:1: 4) & (1:19:19) & (1:4: 1) \\
(1:1: 22) & (1:8: 8) & (1:22:1) \\
(1:3: 12) & (1:23:19) & (1:4: 17) \\
(1:7: 8) & (1:22:4) & (1:19:18) \\
(1:7: 22) & (1:8: 6) & (1:21:18) \\
(5:1: 2) & (1:15:13) & (1:2: 5) \\
(5:1: 23) & (1:10:12) & (1:23:5).
\end{array}
\]

Its automorphism group has order $3$ and is generated by a rotation of the coordinate axes: $\Aut(\arcone) = \left<\rho\right>$ where
\[
\rho = \left<v\right>\mapsto\left<\begin{pmatrix}0 & 1 & 0\\0 & 0 & 1\\1 & 0 & 0\end{pmatrix}v\right>
\]

The automorphism group partitions $\arcone$ into $7$ orbits, each of size $3$.
In the list of points above each row consists of a single orbit.

The $21$ points are contained in $21$ different point neighbor classes.
The complement of $\phi(\arcone)$ in $\PG(2,\F_5)$ has size $10$ and consists of the projective triangle 
\[
\bigcup_{a\in\F_5}\{{(0:1:-a^2)},{(1:-a^2:0)},{(-a^2:0:1)}\}
\]
together with its center point $(1:1:1)$.

\subsection{A $(22,2)$-arc in $\IS_5$}
A $(22,2)$-arc $\arctwo$ in $\PHG(2,\IS_5)$ is given by the points

\[
\begin{array}{ccc}
(1:X+1:4X) & (4X:1:X+1) & (1:4X:4X+1) \\
(1:4X+1:4X) & (4X:1:4X+1) & (1:4X:X+1)\\
(1:X+1:3X+4) & (1:2X+4:X+4) & (1:4X+4:4X+1) \\
(1:4X+1:4X+4) & (1:X+4:2X+4) & (1:3X+4:X+1) \\
\hline
(1:3X+2:3X+2) & (1:3X+3:1) & (1:1:3X+3) \\
(1:2X+3:4X+2) & (1:4X+3:3X+4) & (1:2X+4:2X+2) \\
(1:4X+2:2X+3) & (1:2X+2:2X+4) & (1:3X+4:4X+3) \\
 & (1:1:1). &
\end{array}
\]

Again, the rotation $\rho$ of the coordinate axes is an automorphism of $\arctwo$, and each row in the point list consists of one orbit under the group action of $\left<\rho\right>$, so there are $7$ orbits of size $3$ and the fixed point $(1:1:1)$.
But in this case the full automorphism group $\Aut(\arctwo)$ is bigger than $\left<\rho\right>$, another automorphism of $\arctwo$ is given by
\[
   \tau = \left<v\right> \mapsto \left<\begin{pmatrix}1 & X & -X\\X & 1 & -X\\ 2X+2 & 2X+2 & -X-1\end{pmatrix}v\right>
\]
$\tau$ has order $2$, $\rho\tau$ has order $5$ and together, $\rho$ and $\tau$ generate $\Aut(\arctwo)$: $\Aut(\arctwo) = \left<\rho,\tau\right> \cong \PSL(2,\F_5) \cong A_5$, where $A_5$ denotes the alternating group on a set of size $5$.

While the ring $\Z_{25}$ has a trivial automorphism group, $\Aut(\IS_5)$ is cyclic of order $4$ and generated by the linear extension of $X\mapsto 2X$. So $\PGL(3,\IS_5) \lneq \PGammaL(3,\IS_5)$.
But $\Aut(\arctwo) < \PGL(3,\IS_5)$, so all the automorphisms of $\arctwo$ are purely linear.
Under the action of $\Aut(\arctwo)$, $\arctwo$ splits into $2$ orbits $O_1$ and $O_2$.
$O_1$ has size $12$ and contains the points above the horizontal line in the list, $O_2$ has size $10$ and contains the points below the horizontal line.

The $12$ points in $O_1$ are contained in $6$ point neighbor classes, each class containing $2$ points of $O_1$.
The image of these $6$ point neighbor classes under $\phi$ is the oval
\[
O = \{(0:1:1),(1:0:1),(1:1:0),(-1:1:1),(1:-1:1):(1:1:-1)\}
\]
in $\PG(2,\F_5)$.
Each pair of points within the same point neighbor class is aligned in the tangent direction of the oval $O$.

The $10$ points in $O_2$ are all in seperate point neighbor classes, and their $\phi$-images are exactly the internal points of the oval $O$.
The $15$ point neighbor classes corresponding to the external points of $O$ are empty.

\section{Computation}
The arcs were found by a fast backtracking search.

One problem are the huge automorphism groups of the projective Hjelmslev planes\footnote{$\#\Aut(\PHG(2,\Z_{25})) = 145312500000$ and $\#\Aut(\PHG(2,\IS_5)) = 581250000000$}, which cause ''the same'' point constellation to appear in billions of isomorphic copies during a naive depth-first search.
On the other hand, a complete elimination of isomorphic copies would be too expensive, so the compromise was to filter out isomorphic copies in the first $7$ levels of the search, and to do an additional isomorphism test for the leaf nodes of the search.
These isomorphism tests and the determination of the automorphism groups were done implementing the \emph{Leiterspiel} \cite{Schmalz-1993-BayMS31:109-143}, see also \cite{Laue-1993-BayMS43:53-96}.

Another bottleneck is the test in the innermost loop of the algorithm if a certain point can be added to the current point set without violating the $2$-arc property.
Here we exploit the fact that $\mathfrak{k}\subseteq \mathcal{P}(\PHG(2,R))$ is a $2$-arc if and only if all the $3$-element subsets of $\mathfrak{k}$ are a $2$-arc:
For any set $S\subset\mathcal{P}(\PHG(2,R))$ and any four points $P_1,P_2,P_3,P_4\in\mathcal{P}(\PHG(2,R))\setminus S$ it holds that $\mathfrak{a} := S\cup \{P_1,P_2,P_3,P_4\}$ is a $2$-arc if and only if $\mathfrak{a}_1 := S\cup\{P_1,P_2,P_3\}$, $\mathfrak{a}_2 := S\cup\{P_1,P_2,P_4\}$, $\mathfrak{a}_3 := S\cup\{P_1,P_3,P_4\}$ and $\mathfrak{a}_4 := S\cup\{P_2,P_3,P_4\}$ are $2$-arcs.
So in each node of the depth-first search we do a local breadth-first search for $3$ levels.
The overhead of this additional search is compensated by the fact that the search depth of the outer depth-first search is reduced by $3$.
Now when the backtrack algorithm goes forward from the search node $S$ to the search node $S\cup\{P_1\}$, the breadth-first search for the node $S\cup\{P_1\}$ can be performed easily only by doing look-ups in the breadth-first data of the node $S$.
This process can be seen as \emph{merging} the $4$ arcs $\mathfrak{a}_1$, $\mathfrak{a}_2$, $\mathfrak{a}_3$ and $\mathfrak{a}_4$ into the bigger arc $\mathfrak{a}$.
The merging technique is a general idea to avoid repeated tests within a backtracking search.
In \cite{Kurz-2006-BayMS76} it was used on pairs of integral point sets, and then in \cite{Koch-2006} on triples of polyominoes.

%For the computation of the automorphism groups, we transformed $\PHG(2,R)$ into a bipartite graph on the vertex set $\mathcal{P}(\PHG(2,R)) \cup \mathcal{L}(\PHG(2,R))$, where two vertices $P\in \mathcal{P}(\PHG(2,R))$ and $L\in \mathcal{L}(\PHG(2,R))$ are connected by an edge if and only if the point $P$ is on the line $L$.
%If we assign the same color to all the vertices in $\mathcal{P}(\PHG(2,R))$ and another color to the vertices in $\mathcal{L}(\PHG(2,R))$, the automorphism group of the resulting colored graph is isomorphic to $\Aut(\PHG(2,R))$. And if we further color the vertices corresponding to the points of an arc $\mathfrak{k}$ in yet another color, the automorphism group of this colored graph is isomorphic to $\Aut(\mathfrak{k})$.
%For the computation of the automorphism groups of colored graphs, the software package \verb+nauty+ \cite{McKay-Nauty} was used.

\section{Conclusion and Future Research}
Up to isomorphism, the given $2$-arcs were the only ones of size $21$ respectively $22$ that showed up in our search.
Since we did not investigate the complete search space, there might still exist other isomorphism types or even bigger $2$-arcs.
But the current situation is quite remarkable: When we started our search, in the tables of the best known $u$-arcs in uniform projective Hjelmslev planes over finite chain rings all the $u$-arcs in planes over Galois rings were at least as large as the $u$-arcs in the Hjelmslev planes over the other rings $R$ with the same parameter $q$, suggesting this being true in general.
The Hjelmslev plane $\PHG(2,\IS_5)$ admitting a $2$-arc of size $22$ on the one hand and the Hjelmslev plane $\PHG(2,\Z_{25})$ with its best known arc of size $21$ on the other hand could be a counterexample to this conjecture.

Of course, the definitive knowledge of the biggest $2$-arcs in $\PHG(2,\Z_{25})$ and $\PHG(2,\IS_5)$ would be great.
We think that it might be computationally feasible to exhaustively search the complete search space by further exploiting the homomorphism $\phi$ of group actions
\[(\PGammaL(3,R),\PHG(2,R)) \rightarrow (\PGL(3,\F_5),\PG(2,\F_5))\]
via the homomorphism principle \cite{Laue-2001}:
In a first step, all the $\PGL(3,\F_5)$-representatives for the images in $\PG(2,\F_5)$ are generated.
At this point we can make use of some restrictions, for example that the empty point neighbor classes form a blocking set in $\PG(2,\F_5)$.
Then for each such image $\bar{\mathfrak{k}}$, we need to exhaustivly search all the $2$-arcs among the preimages in $\phi^{-1}(\bar{\mathfrak{k}})$ up to $\PGammaL(3,R)$-isomorphism.
In fact it is enough to consider the preimage of the $\PGL(3,\F_5)$-stablilizer of $\bar{\mathfrak{k}}$ as the operating group.
Usually this group is much smaller than the full group of all collineations, a fact that benefits canonization and isomorphism tests.

Furthermore, the question arises if a generalization of the arc $\arcone$ or $\arctwo$ to uniform Hjelmslev planes over chain rings of higher order is possible.
We hope that our analysis of the structure could be a first step into this direction.

\end{document}